\documentclass[12pt,a4paper]{article}      

\usepackage{graphicx}
\usepackage{enumerate}
\usepackage{amssymb, amsmath, amsthm}

\newtheorem{satz}{Theorem}[section]
\newtheorem{lem}[satz]{Lemma} 
\newtheorem{proposition}[satz]{Proposition} 
 
\newtheorem{cor}[satz]{Corollary}

\usepackage[T1]{fontenc}
\usepackage{lmodern}

\newcommand{\N}{\ensuremath{{\mathbb N}}}

\newcommand{\R}{\ensuremath{{\mathbb R}}}
\newcommand{\e}{\varepsilon}
\newcommand{\AveP}{\underset{\pi}{\mbox{Ave}}}
\newcommand{\AvePP}{\underset{\pi,\sigma}{\mbox{Ave}}}
\newcommand{\AvePPP}{\underset{\pi,\sigma,\eta}{\mbox{Ave}}}
\newcommand{\AvePs}{\underset{\sigma}{\mbox{Ave}}}

\newcommand{\abs}[1]{\left\lvert#1 \right\rvert}
\newcommand{\norm}[1]{\left \lVert#1 \right\rVert}

\begin{document}

\title{Combinatorial Inequalities and Subspaces of $L_1$}
\author{Joscha Prochno \thanks{Part of this paper is part of the doctoral thesis of the first named author  (see \cite{key-6}), supervised by the second named author} \and Carsten Sch\"utt }
\date{\today}
\maketitle

\begin{abstract} 
Let  $M$ and $N$ be Orlicz functions. We establish some combinatorial inequalities and show that the product spaces $\ell^n_M(\ell_N^{n})$ are uniformly isomorphic to subspaces of $L_1$
if $M$ and $N$ are ``separated'' by a function $t^{r}$, $1<r<2$. 
\end{abstract}


\section{Introduction}\label{intro}
The structure and variety of subspaces of $L_1$ is very rich. Over the years, there was put tremendous effort in characterizing subspaces
of $L_{1}$. Although there are a number of sophisticated criteria
at hand now, it might turn out to be nontrivial to decide for
a specific Banach space whether it is isomorphic to a subspace of
$L_{1}$.
\par
Using the theorem of de Finetti
it had been shown in \cite{key-1} that every Orlicz space with a $2$-concave Orlicz function embeds into $L_1$. Consequently, all spaces whose norms are averages of $2$-concave Orlicz norms embed into $L_1$. In fact, this characterizes all subspaces of $L_1$ with a symmetric basis. The corresponding finite-dimensional version of this result was proved in \cite{KS1}, using combinatorial and probabilistic tools.
\par
Although this characterization gives a complete picture of which 
spaces with a symmetric basis embed into $L_{1}$ it might not be 
easy to apply. This becomes apparent when one considers 
Lorentz spaces \cite{Sch}.
\par
Here we study matrix subspaces of $L_1$, {\it i.e.}, spaces $E(F)$ where $E$ and $F$ have a $1$-symmetric basis $(e_i)_{i=1}^n$ and $(f_j)_{j=1}^n$, and where for all matrices $(x_{ij})_{i,j}$
  $$
    \norm{(x_{ij})_{i,j}}_{E(F)}= \norm{\sum_{i=1}^n \norm{\sum_{j=1}^n x_{ij}f_j}_F e_i}_E.
  $$
  
Our main result is the following:

\begin{satz}\label{MatrixEmbed}
Let $1<p<r<2$ and $M$ and $N$ Orlicz functions with
$\frac{M(t)}{t^{p}}$ decreasing and $\frac{N(t)}{t^{r}}$ increasing and $\frac{N(t)}{t^{2}}$
decreasing. Then there is a constant $C>0$ such that for
all $n\in\mathbb N$ there is a subspace $E$ of
$L_{1}$ with $\operatorname{dim}(E)=n^{2}$
and
$$
d(E,\ell_{M}^{n}(\ell_{N}^{n}))\leq C.
$$
\end{satz}
\vskip 5mm

Here, $d$ denotes the Banach-Mazur distance.  To prove Theorem \ref{MatrixEmbed} we first show that $\ell_{M}^n(\ell_r^n)$ is isomorphic to a subspace of $L_1$. To do this, we develop some technical combinatorial results related to Orlicz norms and use techniques first developed in \cite{KS1} and \cite{KS2}. These combinatorial inequalities, used to embed finite-dimensional Banach spaces into $L_1$, are interesting in themselves. Using the results of  Bretagnolle and  Dacunha-Castelle from \cite{key-1}, {\it i.e.}, $\ell_{N}$ is a subspace of $L_{r}$ if and only if $\frac{N(t)}{t^{r}}$ increasing and $\frac{N(t)}{t^{2}}$ decreasing, we obtain our main result. 
\par
In some sense the conditions that $\frac{M(t)}{t^p}$ is decreasing, $\frac{N(t)}{t^r}$ is increasing and $\frac{N(t)}{t^2}$ is decreasing, are strict. This is a consequence of a result from \cite{KS2} (Corollary 3.3). Kwapie\'n and Sch\"utt proved that
  $$
    \frac{1}{5\sqrt{2}} \norm{Id} \leq d\left(E(F),G\right),
  $$
where $Id\in L(E,F)$ is the natural identity map, i.e. $Id(\sum_{i=1}^n a_ie_i)=\sum_{j=1}^n a_jf_j$, and $E$, $F$ are $n$-dimensional spaces with a 
$1$-symmetric and $1$-unconditional basis respectively. For $1\leq p<r \leq 2$ they obtain that for any $n^2$-dimensional subspace $G$ of $L_1$ 
  $$
    d\left(\ell^n_r(\ell_p^n), G\right) \geq \frac{1}{5\sqrt{2}}n^{1/p-1/r}
  $$   
holds. Therefore, the conditions are strict.
\par
The technical difficulties that occur are that in general Orlicz functions are not homogeneous for some $p$, {\it i.e.}, $M(\lambda t)\neq \lambda^p M(t)$.
\par
Furthermore, since our results are of a very technical nature in many places, we tried to make this paper as self contained as possible and therefore easily accessible. 

\section{Preliminaries and combinatorial inequalities}

A convex function $M:[0,\infty) \to [0,\infty)$ with $M(0)=0$ and $M(t)>0$ for $t>0$ is called an Orlicz function. We define the Orlicz space $\ell^n_M$ to be $\R^n$ equipped with the norm
  $$
    \norm{x}_M = \inf\left\{ \rho>0 \Big| \sum_{i=1}^n M\left(\frac{\abs{x_i}}{\rho}\right) \leq 1  \right\}.
  $$
Given an Orlicz function $M$, we define its dual function $M^*$ by the Legendre transform, {\it i.e.},
  $$
    M^*(x) = \sup_{t\in[0,\infty)}(xt-M(t)).
  $$
We have for all Orlicz functions $M$ and $0\leq t <\infty$
  \begin{equation}\label{EQU M hoch -1 M stern hoch -1}
    t \leq M^{-1}(t) M^{*-1}(t) \leq 2t,
  \end{equation}
A proof can be found in \cite{key-2} (Formula 2.10, page 13). 
We say that two Orlicz functions $M$ and $N$ are equivalent
if there are positive constants $a$ and $b$ such that for all
$t\geq0$
$$
aN^{-1}(t)\leq M^{-1}(t)\leq bN^{-1}(t).
$$
If two Orlicz functions are equivalent so are their norms.
\par
Let $X$ and $Y$ be isomorphic Banach spaces. We say that they are
$C$-isomorphic if there is an isomorphism $I:X\rightarrow Y$ with
$\|I\|\|I^{-1}\|\leq C$.
We define the Banach-Mazur distance of $X$ and $Y$ by
    $$
      d(X,Y) = \inf\left\{ \norm{T}\norm{T^{-1}} \Big| ~ T\in L(X,Y) ~ \hbox{isomorphism} \right\}.
    $$
 Let $(X_n)_n$ be a sequence of $n$-dimensional normed spaces and let $Z$ also be a normed space. If there exists a constant $C>0$, such that for all $n\in\N$ there exists a normed space $Y_n \leq Z$ with $\dim(Y_n)=n$ and $d(X_n,Y_n)\leq C$, then we say that $(X_n)_n$ embeds uniformly into $Z$ or in short: $X_n$ embeds into $Z$. For a detailed introduction to the concept of Banach-Mazur distances, see for example \cite{key-5}.  

\par
We need the following two results by Kwapie\'n and Sch\"utt from \cite{KS1, KS2}.

\begin{lem} [\cite{KS1} Lemma 2.1] \label{LEM Kwapien-Schuett 1}
  Let $n,m\in\N$ with $n \leq m$ and let $y\in\R^m$ with $y_1\geq y_2 \geq \ldots \geq y_m >0$. Furthermore, let $M$ be an Orlicz function such 
  that for all $k=1,\ldots,m$
    $$
      M^*\left( \sum_{i=1}^k y_i \right)= \frac{k}{m}.
    $$
  We define $\norm{\cdot}_y$ by
    $$
      \norm{x}_y = \max_{\sum_{i=1}^nk_i=m} \sum_{i=1}^n \left( \sum_{j=1}^{k_i}y_j \right) |x_i|.
    $$
  Then, for all $x\in\R^n$,
    $$
      \frac{1}{2} \norm{x}_y \leq \norm{x}_M \leq 2 \norm{x}_y.
    $$    
\end{lem}

\vskip 5mm

\begin{lem} [\cite{KS2} Lemma 2.5] \label{LEM Kwapien-Schuett 2}
  Let $M$ be an Orlicz function. Then, for all $x\in\R^n$,
    \begin{eqnarray*}
      && \frac{1}{2}\left(\frac{1}{2}-\frac{1}{n-1}\right)\norm{x}_M \\
      & \leq & \frac{1}{n!}\sum_{\pi} \max_{1\leq i \leq n} \left| x_i\cdot n \cdot \left(M^{*-1}\left(\frac{\pi(i)}{n}\right)
      - M^{*-1}\left(\frac{\pi(i)-1}{n}\right) \right) \right| ~ \leq ~ 2 \norm{x}_M .
    \end{eqnarray*}
\end{lem}

\vskip 5mm

\begin{lem}[\cite{KS2} Corollary 1.7] \label{3MatrixMittel}
For all $n\in\mathbb N$ and all nonnegative numbers
$B(i,k,\ell)$, $1\leq i,k,\ell\leq n$,
$$
\frac{1}{16 n^{2}}
\sum_{\alpha=1}^{n^{2}}s(\alpha)
\leq\frac{1}{(n!)^{2}}
\sum_{\pi,\sigma}\max_{1\leq i\leq n} B(i,\pi(i),\sigma(i))
\leq\frac{4}{n^{2}}
\sum_{\alpha=1}^{n^{2}}s(\alpha),
$$
where $s(1),\dots,s(n^{3})$ is the decreasing rearrangement
of the numbers $B(i,k,\ell)$, $1\leq i,k,\ell\leq n$.\end{lem}
\vskip 3mm

From Lemma \ref{LEM Kwapien-Schuett 1} and Lemma \ref{LEM Kwapien-Schuett 2} we obtain the following result.

\begin{lem} \label{LEM M-Norm liefert N-Norm}
  Let $a\in\R^n$, $a_1 \geq \ldots \geq a_n > 0$ and let $M$ be an Orlicz function. Furthermore, let $N$ be an Orlicz function such that for 
  the dual function $N^*$ and all $\ell=1,\ldots,n^2$
    $$
      N^{*-1}\left( \frac{\ell}{n^2} \right) = \frac{1}{n^2} \sum_{k=1}^{\ell} s(k)
    $$
  holds, where $s(1),\ldots,s(n^2)$ is the decreasing rearrangement of
    $$
      a_i\cdot n\cdot \left(M^{*-1}\left(\frac{j}{n}\right)-M^{*-1}\left(\frac{j-1}{n}\right) \right),~i,j=1,\ldots,n.
    $$
  Then, for all $x\in\R^n$,
    $$
      c \norm{x}_N \leq \frac{1}{n!} \sum_{\pi} \norm{(x_ia_{\pi(i)})_{i=1}^n}_M \leq 2\norm{x}_N,
    $$
  where $c>0$ is an absolute constant. \\  
  Furthermore, one can choose $N$, such that $N^{*-1}$ is an affine function between the values $\frac{\ell}{n^2}$, $\ell=1,\ldots,n^2$.     
\end{lem}

\noindent
{\bf Proof.}
  From Lemma \ref{LEM Kwapien-Schuett 2}, we know
    $$
      c \norm{x}_M \leq 
      \frac{1}{n!} \sum_{\sigma} \max_{1\leq i \leq n} \left| x_i\cdot n \cdot \left(M^{*-1}\left(\frac{\sigma(i)}{n}\right)
      - M^{*-1}\left(\frac{\sigma(i)-1}{n}\right) \right) \right| \leq  2 \norm{x}_M.
    $$
  Thus
    \begin{eqnarray*}
      && c \frac{1}{n!}\sum_{\pi} \norm{(x_i a_{\pi(i)})_{i=1}^n}_M\\
      & \leq & \frac{1}{n!^2} \sum_{\sigma,\pi} \max_{1\leq i \leq n} \left| x_ia_{\pi(i)}\cdot n \cdot \left(M^{*-1}\left(\frac{\sigma(i)}{n}\right)
      - M^{*-1}\left(\frac{\sigma(i)-1}{n}\right) \right) \right| \\
      & \leq & 2 \frac{1}{n!} \sum_{\pi} \norm{(x_ia_{\pi(i)})_{i=1}^n}_M.
    \end{eqnarray*} 
  Applying Lemma \ref{LEM Kwapien-Schuett 1} and \ref{3MatrixMittel} yields the desired result.   
$\square$

\vskip 5mm

Now we are able to develop the combinatorial ingredients that we need to prove Theorem \ref{THM Joscha u Carsten 2}.  These results are extensions of results that had been shown in \cite{KS1}, respectively \cite{KS2}.

\begin{lem}\label{KombMittelLr}
(i) Let $1< r<\infty$ and 
$a_{1}\geq a_{2}\geq\cdots\geq a_{n}>0$.
Then there exists an Orlicz function
$N$ such that for the dual function $N^{*}$ and all 
$\ell=1,\dots,n$
\begin{equation}\label{KombMittelLr6}
N^{*-1}\left(\frac{\ell}{n}\right)
\leq C_r
\left(\left(\frac{1}{n}
\sum_{i=1}^{\ell}a_{i}\right)
+\left(\frac{\ell}{n}\right)^{\frac{1}{r^{*}}}
\left(\frac{1}{n}\sum_{i=\ell+1}^{n}|a_{i}|^{r}
\right)^{\frac{1}{r}}\right)
\leq 8N^{*-1}\left(\frac{\ell}{n}\right)
\end{equation}
\begin{equation}\label{KombMittelLr7}
N^{*-1}\left(\frac{\ell}{n^{2}}\right)
\leq
C_r
\frac{1}{n}\left(\frac{\ell}{n}\right)^{\frac{1}{r^{*}}}
\left(\sum_{i=1}^{\ell}|a_{i}|^{r}\right)^{\frac{1}{r}}
\leq
2N^{*-1}\left(\frac{\ell}{n^{2}}\right),
\end{equation}
where $C_r=r^{\frac{1}{r}}\left(r^{*}\right)^{\frac{1}{r^{*}}}$. Furthermore, for all $x\in\mathbb R^{n}$ 
$$
c\|x\|_{N}
\leq \frac{1}{n!}\sum_{\pi}
\left(\sum_{i=1}^{n}|x_{i}a_{\pi(i)}|^{r}\right)^{\frac{1}{r}}
\leq 2\|x\|_{N}.
$$
\newline
(ii) Let $1< r<\infty$. There exists $n_{0}\in\mathbb N$ such that for all $n\geq n_{0}$, all
$a_{1}\geq \ldots \geq a_{n}>0$ and all Orlicz functions $\bar N$, where for all $\ell=1,\dots,n$
\begin{equation}\label{KombMittelLr71}
\bar N^{*-1}\left(\frac{\ell}{n}\right)
\leq C_r
\left(\left(\frac{1}{n}
\sum_{i=1}^{\ell}a_{i}\right)
+\left(\frac{\ell}{n}\right)^{\frac{1}{r^{*}}}
\left(\frac{1}{n}\sum_{i=\ell+1}^{n}|a_{i}|^{r}
\right)^{\frac{1}{r}}\right)
\leq 8\bar N^{*-1}\left(\frac{\ell}{n}\right)
\end{equation}
holds and which are affine on the intervals 
$[\frac{\ell}{n},\frac{\ell+1}{n}]$, $\ell=0,\dots,n-1$,
we have for all $x\in\mathbb R^{n}$
$$
a_{r}\|x\|_{\bar N}
\leq \frac{1}{n!}\sum_{\pi}
\left(\sum_{i=1}^{n}|x_{i}a_{\pi(i)}|^{r}\right)^{\frac{1}{r}}
\leq b_{r}\|x\|_{\bar N},
$$
where $a_{r}$ and $b_{r}$ just depend on $r$ and $C_r$ as in $(i)$.
\end{lem}
\vskip 3mm

By part (i) there is indeed an Orlicz function as it is specified in (ii):
The Orlicz function as given in (i) can be modified so that it is affine 
on the intervals $[\frac{\ell}{n},\frac{\ell+1}{n}]$, $\ell=0,\dots,n-1$.
\vskip 3mm

\noindent
{\bf Proof.}
(i) From Lemma \ref{LEM M-Norm liefert N-Norm} we obtain
$$
c_{}\|x\|_{N}
\leq\frac{1}{n!}\sum_{\pi}\|(x(i)a_{\pi(i)})_{i=1}^{n}\|_{M}
\leq 2\|x\|_{N},
$$
where $M(t)=t^{r}$,
$$
N^{*}\left(
\frac{1}{n^{2}}\sum_{k=1}^{\ell}s(k)
\right)=\frac{\ell}{n^{2}},~
\ell=1,\dots,n^{2}
$$
and
$s(1),\dots,s(n^{2})$ is the decreasing rearrangement of the numbers
$$
a_{i}\cdot n\cdot
\left(M^{*-1}(\tfrac{j}{n})-M^{*-1}(\tfrac{j-1}{n})\right),~
1\leq i,j\leq n.
$$
Obviously $M^{*}(s)=(\frac{1}{r})^{\frac{1}{r}}
(\frac{1}{r^{*}})^{\frac{1}{r^{*}}}s^{r^{*}}$ and
$
M^{*-1}(t)
=r^{\frac{1}{r}}\left(r^{*}\right)^{\frac{1}{r^{*}}}
t^{\frac{1}{r^{*}}}
$. We choose $C_r:=r^{\frac{1}{r}}\left(r^{*}\right)^{\frac{1}{r^{*}}}$. For all $\ell\leq n^2$ we have
\begin{eqnarray*}
  \frac{1}{n^{2}}\sum_{k=1}^{\ell}s(k)
  &=&\max_{\begin{array}{c}
  \sum_{i=1}^{n}\ell_{i}=\ell   \\
  \ell_{i}\leq n
  \end{array}}
  \frac{1}{n^{2}}\sum_{i=1}^{n}a_{i}
  \sum_{j=1}^{\ell_{i}}n\cdot
  \left(M^{*-1}(\tfrac{j}{n})-M^{*-1}(\tfrac{j-1}{n})\right)
  \\
  &=&\max_{\begin{array}{c}
  \sum_{i=1}^{n}\ell_{i}=\ell   \\
  \ell_{i}\leq n
  \end{array}}
  \frac{1}{n^{2}}\sum_{i=1}^{n}a_{i}
  n\cdot
  M^{*-1}(\tfrac{\ell_{i}}{n})
  \\
  &=&\max_{\begin{array}{c}
  \sum_{i=1}^{n}\ell_{i}=\ell   \\
  \ell_{i}\leq n
  \end{array}}
  C_r \frac{1}{n^{}}\sum_{i=1}^{n}a_{i}
  |\tfrac{\ell_{i}}{n}|^{\frac{1}{r^{*}}}.
\end{eqnarray*}
Thus, for all $\ell \leq n^2$
\begin{equation}\label{KombMittelLr1}
\frac{1}{n^{2}}\sum_{k=1}^{\ell}s(k)
=\max_{\begin{array}{c}
\sum_{i=1}^{n}\ell_{i}=\ell   \\
\ell_{i}\leq n
\end{array}}
C_r \frac{1}{n^{}}\sum_{i=1}^{n}a_{i}
|\tfrac{\ell_{i}}{n}|^{\frac{1}{r^{*}}}.
\end{equation}
We show the right hand side inequality  (\ref{KombMittelLr6}). 
We consider the case $\ell=m\cdot n$, $1\leq m \leq n$. Then
  $$
    N^{*-1}\left(\frac{m}{n}\right) = \frac{1}{n^2} \sum_{k=1}^{n\cdot m}s(k),
  $$
and by (\ref{KombMittelLr1})
  $$
    N^{*-1}\left(\frac{m}{n}\right) = C_r \max_{\begin{array}{c} \sum_{i=1}^{n}\ell_{i}=m\cdot n   
    \\ \ell_{i}\leq n \end{array}} \frac{1}{n^{}}\sum_{i=1}^{n}a_{i}|\tfrac{\ell_{i}}{n}|^{\frac{1}{r^{*}}}.
  $$ 
  For $m=1$ we obtain from Lemma \ref{LEM Kwapien-Schuett 1} that $ N^{*-1}\left(\frac{1}{n}\right)$
  is of the order
  $\|a\|_{r}$ . Now we consider $m\geq 2$.
We choose $\ell_{1}=\ldots=\ell_{m}=n$ and 
$\ell_{m+1}=\ldots=\ell_{n}=0$ and obtain
  \begin{equation}\label{KombMittelLr17}
    N^{*-1}\left(\frac{m}{n}\right)
    \geq
    C_r \frac{1}{n^{}}\sum_{i=1}^{m}a_{i}.
  \end{equation}
We consider
  $$
    y_j=M^{*-1}\left(\frac{j}{n\cdot m}\right) - M^{*-1}\left( \frac{j-1}{n\cdot m} \right),~~1\leq j \leq n\cdot m.
  $$
From Lemma \ref{LEM Kwapien-Schuett 1} we get
  \begin{eqnarray*}
    \frac{1}{2}\norm{a}_r & \leq & \norm{a}_y 
    ~ = ~ \max_{\sum_{i=1}^{n}\ell_{i}=m\cdot n   }
    \sum_{i=1}^{n}a_{i} \left(\sum_{j=1}^{l_i} y_j\right) \\
    & = & C_r \max_{\sum_{i=1}^{n}\ell_{i}=m\cdot n}
    \sum_{i=1}^{n}a_{i}
    |\tfrac{\ell_{i}}{nm}|^{\frac{1}{r^{*}}}. \\
  \end{eqnarray*}  
This holds if and only if 
  $$
    \frac{1}{2}m^{1/r^*}\norm{a}_r 
    \leq  C_r \max_{\sum_{i=1}^{n}\ell_{i}=m\cdot n}
    \sum_{i=1}^{n}a_{i}
    |\tfrac{\ell_{i}}{n}|^{\frac{1}{r^{*}}}.
  $$
The inequality also holds for the modified vector
$\tilde a$ with $\tilde a_{1}=\ldots=\tilde a_{m}=a_{m}$
and $\tilde a_{i}=a_{i}$ for $i=m+1,\dots,n$, i.e.
  \begin{equation} \label{EQU Normabschaetzung fuer Vektor tilde a}
    \frac{1}{2}m^{\frac{1}{r^{*}}}\|\tilde a\|_{r}
    \leq C_r \max_{\sum_{i=1}^{n}\ell_{i}=m\cdot n}
    \sum_{i=1}^{n}\tilde a_{i}
    |\tfrac{\ell_{i}}{n}|^{\frac{1}{r^{*}}}.
  \end{equation}
We show that w.l.o.g. $\ell_i \leq n$, $i\leq n$. It holds $\ell_1 \geq \ldots \geq \ell_m \geq \ldots \geq \ell_n \geq 0$. Obviously, we have $\ell_i\leq n$ for all $i=m,\ldots,n$. Otherwise we would have $\sum_{i=1}^m \ell_i>n\cdot m$, which cannot occur. Therefore, it suffices to show that we can choose $\ell_1,\ldots, \ell_m \leq n$. To do this, we construct $\tilde \ell_i$, $i\leq m$ such that $\tilde \ell_i \leq n$, and such that the maximum in (\ref{EQU Normabschaetzung fuer Vektor tilde a}) is attained up to an absolute constant (we take $\tilde \ell_i=\ell_i$ for $i=m+1,\ldots,n$). Now, let $\ell_1,\ldots,\ell_n$ such that the maximum in (\ref{EQU Normabschaetzung fuer Vektor tilde a}) is attained. Then we define for $i\leq m$
  $$
    \tilde \ell_i := \left\lfloor \frac{1}{m} \sum_{j=1}^m \ell_j \right\rfloor.
  $$  
  ($\lfloor x \rfloor$ is the biggest integer smaller than $x$.)
We may assume 
  $$
    \left\lfloor \frac{1}{m} \sum_{j=1}^m \ell_j \right\rfloor \geq 2,
  $$  
because otherwise from $\left\lfloor \frac{1}{m} \sum_{j=1}^m \ell_j \right\rfloor < 2$ we obtain immediately that $\ell_{m+1},\ldots,\ell_n\leq 1$, and therefore
$n\cdot m=\sum_{j=1}^{n}\ell_{j}<2m+(n-m)=n+m$. 
Since $m\geq2$ and we may assume that $n\geq 3$ we get
a contradiction. 
Hence, we get for all $i\leq m$
  $$
    \tilde \ell_i \geq \frac{1}{m} \sum_{j=1}^m \ell_j - 1 \geq \frac{1}{2} \frac{1}{m} \sum_{j=1}^m \ell_j. 
  $$  
Now we have
  $$
    \sum_{i=1}^m \tilde a_i | \tfrac{\tilde \ell_i}{n} |^{1/r^*} \geq 
    \sum_{i=1}^m a_m \frac{1}{2^{1/r^*}} \left( \frac{1}{n} \frac{1}{m} \sum_{j=1}^m \ell_j \right)^{1/r^*}=
    a_m m^{1/r} \frac{1}{2^{1/r^*}} \left( \frac{1}{n} \sum_{j=1}^m \ell_j \right)^{1/r^*}.
  $$
From H\"older's inequality we get  
  $$
    \sum_{i=1}^{m}\tilde a_{i}
    |\tfrac{\ell_{i}}{n}|^{\frac{1}{r^{*}}}
    =a_{m}\sum_{i=1}^{m}
    |\tfrac{\ell_{i}}{n}|^{\frac{1}{r^{*}}}
    \leq a_{m}m^{\frac{1}{r}}\left(\frac{1}{n}\sum_{i=1}^{m}\ell_{i}
    \right)^{\frac{1}{r^{*}}}.
  $$
Thus
  $$
    \sum_{i=1}^m \tilde a_i | \tfrac{\tilde \ell_i}{n} |^{1/r^*} \geq \frac{1}{2^{1/r^*}} \sum_{i=1}^{m}\tilde a_{i}
    |\tfrac{\ell_{i}}{n}|^{\frac{1}{r^{*}}},
  $$
and therefore
  $$
    \sum_{i=1}^n \tilde a_i | \tfrac{\tilde \ell_i}{n} |^{1/r^*} \geq \frac{1}{2^{1/r^*}} \sum_{i=1}^{n}\tilde a_{i}
    |\tfrac{\ell_{i}}{n}|^{\frac{1}{r^{*}}}.
  $$
Inequality (\ref{EQU Normabschaetzung fuer Vektor tilde a}) gives us
  $$
    C_r \sum_{i=1}^n \tilde a_i | \tfrac{\tilde \ell_i}{n} |^{1/r^*} 
    \geq \frac{1}{2^{1/r^*}}\frac{1}{2}m^{1/r^*} \norm{\tilde a}_r.
  $$
So we have
  $$
    \frac{1}{2^{1/r^*}}\frac{1}{2}m^{1/r^*} \norm{\tilde a}_r
    \leq C_r \max_{\begin{array}{c}
    \sum_{i=1}^{n}\tilde \ell_{i}=m\cdot n   \\
    \tilde \ell_{i}\leq n 
    \end{array}} \sum_{i=1}^n \tilde a_i | \tfrac{\tilde \ell_i}{n} |^{1/r^*}
    = n N^{*-1}\left(\frac{m}{n}\right),
  $$
{\it i.e.}, 
  $$
    N^{*-1}\left(\frac{m}{n}\right) \geq \frac{1}{2^{1/r^*}}\frac{1}{2}\frac{m^{1/r^*}}{n} \norm{\tilde a}_r,
  $$  
and because of
  $$
    \left( \sum_{i=m+1}^n |a_i|^r \right)^{1/r} \leq \|\tilde a\|_{r},
  $$
  and (\ref{KombMittelLr17})
we obtain the right hand side of inequality (\ref{KombMittelLr6}). 
\par
Now we give the estimate on the left hand side of (\ref{KombMittelLr6}).
By (\ref{KombMittelLr1}) for a suitable choice of $\ell_{i}$ 
$$
\frac{1}{n^{2}}\sum_{k=1}^{m\cdot n}s(k)
=C_r
\frac{1}{n^{}}\sum_{i=1}^{n}a_{i}
|\tfrac{\ell_{i}}{n}|^{\frac{1}{r^{*}}}.
$$
Since $\ell_{i}\leq n$, we obtain
$$
\frac{1}{n^{2}}\sum_{k=1}^{m\cdot n}s(k)
\leq C_r
n^{-\frac{1}{r^{*}}-1}
\left(\sum_{i=1}^{m}a_{i}n^{\frac{1}{r^{*}}}
+\sum_{i=m+1}^{n}a_{i}\ell_{i}^{\frac{1}{r^{*}}}\right).  
$$
H\"older's inequality implies
$$
\frac{1}{n^{2}}\sum_{k=1}^{m\cdot n}s(k)
\leq C_r
\frac{1}{n}\left\{\sum_{i=1}^{m}a_{i}
+n^{-\frac{1}{r^{*}}}
\left(\sum_{i=m+1}^{n}|a_{i}|^{r}\right)^{\frac{1}{r}}
\left(\sum_{i=m+1}^{n}\ell_{i}
\right)^{\frac{1}{r^{*}}}\right\}  
$$
and with $\sum_{i=1}^{n}\ell_{i}=m\cdot n$
$$
\frac{1}{n^{2}}\sum_{k=1}^{m\cdot n}s(k)
\leq C_r
\frac{1}{n}\left(\sum_{i=1}^{m}a_{i}
+m^{\frac{1}{r^{*}}}
\left(\sum_{i=m+1}^{n}|a_{i}|^{r}\right)^{\frac{1}{r}}
\right).
$$
Therefore, we obtain the inequality on the left hand side of (\ref{KombMittelLr6}), i.e.
$$
N^{*-1}\left(\frac{m}{n}\right)
\leq C_r \left(\frac{1}{n}
\sum_{i=1}^{m}a_{i}\right)
+\left(\frac{m}{n}\right)^{\frac{1}{r^{*}}}
\left(\frac{1}{n}\sum_{i=m+1}^{n}|a_{i}|^{r}
\right)^{\frac{1}{r}}.
$$
Now we prove inequality (\ref{KombMittelLr7}), i.e.
$$
N^{*-1}\left(\frac{\ell}{n^{2}}\right)
\leq
C_r
\frac{1}{n}\left(\frac{\ell}{n}\right)^{\frac{1}{r^{*}}}
\left(\sum_{i=1}^{\ell}|a_{i}|^{r}\right)^{\frac{1}{r}}
\leq
2N^{*-1}\left(\frac{\ell}{n^{2}}\right).
$$
Because $m\leq n$, from (\ref{KombMittelLr1}) we get 
$$
\frac{1}{n^{2}}\sum_{k=1}^{m}s(k)
=\max_{
\sum_{i=1}^{n}\ell_{i}=m}
C_r
\frac{1}{n^{}}\sum_{i=1}^{n}a_{i}
|\tfrac{\ell_{i}}{n}|^{\frac{1}{r^{*}}}
=
C_r
\frac{1}{n^{}}\left(\frac{m}{n}\right)^{\frac{1}{r^{*}}}
\max_{\sum_{i=1}^{m}\ell_{i}=m}
\sum_{i=1}^{n}a_{i}
|\tfrac{\ell_{i}}{m}|^{\frac{1}{r^{*}}}.
$$
For $m=1,\dots,n$, using H\"older's inequality, we get
the left hand side inequality of (\ref{KombMittelLr7})
$$
N^{*-1}\left(\frac{m}{n^{2}}\right)
=\frac{1}{n^{2}}\sum_{k=1}^{m}s(k)
\leq
C_r
\frac{1}{n^{}}\left(\frac{m}{n}\right)^{\frac{1}{r^{*}}}
\left(\sum_{i=1}^{n}|a_{i}|^{r}\right)^{\frac{1}{r}}.
$$
From Lemma \ref{LEM Kwapien-Schuett 1} we obtain for
$m=1,\dots,n$ the right hand side inequality
of  (\ref{KombMittelLr7})
$$
\frac{1}{n^{}}\left(\frac{m}{n}\right)^{\frac{1}{r^{*}}}
\left(\sum_{i=1}^{n}|a_{i}|^{r}
\right)^{\frac{1}{r}}
\leq \frac{2}{C_r}
N^{*-1}\left(\frac{m}{n^{2}}\right).
$$
\par
(ii) Let $N$ be an Orlicz function as given by part (i).
We show that for all $t$ with $\frac{1}{4^{r^{*}}n}\leq t\leq 1$
\begin{equation}\label{KombMittelLr74}
\frac{1}{8\cdot 4^{r^{*}}}N^{*-1}(t)
\leq \bar N^{*-1}(t)
\leq 32\cdot 4^{r^{*}}N^{*-1}(t).
\end{equation}
From this it follows that for all $x$ 
\begin{equation}\label{KombMittelLr73}
\frac{1}{32\cdot4^{r^{*}}}\|x\|_{N}\leq\|x\|_{\bar N}\leq (48\cdot 4^{r^{*}}+16)\|x\|_{N}.
\end{equation}
We show the estimates first for $\frac{1}{n}\leq t\leq1$.
The range $\frac{1}{4^{r^{*}}n}\leq t\leq\frac{1}{n}$
will be considered after that.
For all $t$ with $\frac{1}{n}\leq t\leq1$, we have
\begin{equation}\label{KombMittelLr4}
\frac{1}{16}N^{*-1}(t)
\leq \bar N^{*-1}(t)
\leq 16N^{*-1}(t).
\end{equation}
We show this.
There exists an $\ell\in\{1,\dots,n-1\}$ such that
$\frac{\ell}{n}\leq t\leq\frac{\ell+1}{n}$. 
By (\ref{KombMittelLr6}) and (\ref{KombMittelLr71})
\begin{eqnarray*}
N^{*-1}(t)
& \leq & N^{*-1}\left(\frac{\ell+1}{n}\right)
\leq C_r \left\{ \left(\frac{1}{n}
\sum_{i=1}^{\ell+1}a_{i}\right)
+\left(\frac{\ell+1}{n}\right)^{\frac{1}{r^{*}}}
\left(\frac{1}{n}\sum_{i=\ell+2}^{n}|a_{i}|^{r}
\right)^{\frac{1}{r}} \right\}   \\
& \leq &
2 C_r \left\{\left(\frac{1}{n}
\sum_{i=1}^{\ell}a_{i}\right)
+\left(\frac{\ell}{n}\right)^{\frac{1}{r^{*}}}
\left(\frac{1}{n}\sum_{i=\ell+1}^{n}|a_{i}|^{r}
\right)^{\frac{1}{r}}\right\} \\  
& \leq & 16\bar N^{*-1}\left(\tfrac{\ell}{n}\right)
~\leq~ 16\bar N^{*-1}(t).
\end{eqnarray*}
The inverse estimate is obtained in the same way.
We show that for all $t$, where
$\frac{1}{4^{r^{*}}n}\leq t\leq \frac{1}{n}$
\begin{equation}\label{KombMittelLr2}
\frac{1}{8\cdot 4^{r^{*}}}N^{*-1}(t)
\leq \bar N^{*-1}(t)
\leq 32\cdot4^{r^{*}}N^{*-1}(t).
\end{equation}
By (\ref{KombMittelLr6}) for $\ell=1$
$$
N^{*-1}(\frac{1}{n})
\leq
C_r
\left(\frac{a_{1}}{n}
+\left(\frac{1}{n}\right)^{\frac{1}{r^{*}}}
\left(\frac{1}{n}\sum_{i=2}^{n}|a_{i}|^{r}
\right)^{\frac{1}{r}}\right)  
$$
For $n$ with $n\geq 2\cdot 4^{r^*}$ we have
$2\cdot4^{r^{*}}[\frac{n}{4^{r^{*}}}]\geq n$.
By H\"older's inequality
\begin{eqnarray*}
N^{*-1}(\frac{1}{n})
&\leq &
2^{\frac{1}{r^{*}}}
C_r
\frac{1}{n}
\left(\sum_{i=1}^{n}|a_{i}|^{r}
\right)^{\frac{1}{r}}
\leq C_{r}2^{\frac{1}{r^{*}}}2\cdot4^{r^{*}}\frac{1}{n}
\left(\sum_{i=1}^{[\frac{n}{4^{r^{*}}}]}|a_{i}|^{r}
\right)^{\frac{1}{r}}   \\
&\leq&  C_{r}2^{\frac{1}{r^{*}}}2\cdot4^{r^{*}}
N^{*-1}\left(\frac{[\frac{n}{4^{r^{*}}}]}{n^{2}}\right)
\leq  C_{r}2^{\frac{1}{r^{*}}}2\cdot4^{r^{*}}
N^{*-1}\left(\frac{1}{4^{r^{*}}n}\right) \\
& \leq &  C_{r}2^{\frac{1}{r^{*}}}2\cdot4^{r^{*}}N^{*-1}(t).
\end{eqnarray*}
Thus, we have for all
$\frac{1}{4^{r^{*}}n}\leq t\leq \frac{1}{n}$
\begin{equation}\label{KombMittelLr41}
N^{*-1}(\frac{1}{n})
\leq  C_{r}2^{\frac{1}{r^{*}}}2\cdot4^{r^{*}}N^{*-1}(t).
\end{equation}
The function $\bar N^{*-1}$ takes the values $\bar N^{*-1}(t)=tn\bar N^{*-1}(\frac{1}{n})$ on the interval
$[0,\frac{1}{n}]$. Hence, for all $t$ with 
$\frac{1}{4^{r^{*}}n}\leq t\leq \frac{1}{n}$
$$
\bar N^{*-1}\left(\frac{1}{4^{r^{*}}n}\right)
\leq\bar  N^{*-1}(t)=tn\bar N^{*-1}\left(\frac{1}{n}\right)
\leq\bar  N^{*-1}\left(\frac{1}{n}\right)
=4^{r^{*}}\bar  N^{*-1}\left(\frac{1}{4^{r^{*}}n}\right).
$$
Thus, we have
\begin{eqnarray*}
N^{*-1}(t)
& \leq & N^{*-1}(\frac{1}{n})
~ \leq ~ C_r \left(\frac{a_{1}}{n}
+\left(\frac{1}{n}\right)^{\frac{1}{r^{*}}}
\left(\frac{1}{n}\sum_{i=2}^{n}|a_{i}|^{r}
\right)^{\frac{1}{r}}\right) \\
& \leq & 8
\bar N^{*-1}\left(\frac{1}{n}\right)
~ \leq ~ 8\cdot 4^{r^{*}}\bar N^{*-1}(t)
\end{eqnarray*}
and
\begin{eqnarray*}
\bar N^{*-1}(t)
& \leq & \bar  N^{*-1}(\frac{1}{n})
\leq C_r \left(\frac{a_{1}}{n}
+\left(\frac{1}{n}\right)^{\frac{1}{r^{*}}}
\left(\frac{1}{n}\sum_{i=2}^{n}|a_{i}|^{r}
\right)^{\frac{1}{r}}\right) \\
& \leq & 8 N^{*-1}\left(\frac{1}{n}\right)
~ \leq ~ C_{r}2^{\frac{1}{r^{*}}}16\cdot4^{r^{*}}N^{*-1}(t).
\end{eqnarray*}
Hence, (\ref{KombMittelLr2}) follows.
Furthermore, we have
\begin{equation}\label{KombMittelLr61}
N^{*-1}\left(\frac{1}{4^{r^{*}}n}\right)
\leq\frac{2}{3}N^{*-1}\left(\frac{1}{n}\right)
\hskip 5mm \hbox{and} \hskip 5mm
\bar N^{*-1}\left(\frac{1}{4^{r^{*}}n}\right)
=\frac{1}{4^{r^{*}}}
\bar N^{*-1}\left(\frac{1}{n}\right).
\end{equation}
The equality is obvious. We show the inequality.
By (\ref{KombMittelLr7}) we get
for $\ell=[\frac{n}{4^{r^{*}}}]$ 
\begin{eqnarray*}
\frac{1}{C_r}
N^{*-1}\left(\frac{1}{4^{r^{*}}n}\right)
&\leq&\frac{1}{C_r}
N^{*-1}\left(\frac{[\frac{n}{4^{r^{*}}}]+1}{n^{2}}\right)  \\
&\leq&
\frac{1}{n}\left(\frac{[\frac{n}{4^{r^{*}}}]+1}{n}\right)^{\frac{1}{r^{*}}}
\left(\sum_{i=1}^{[\frac{n}{4^{r^{*}}}]+1}|a_{i}|^{r}
\right)^{\frac{1}{r}}    .
\end{eqnarray*}
By (\ref{KombMittelLr7}), we get for $\ell=n$ and 
for sufficiently big $n$
$$
\frac{1}{C_r}
N^{*-1}\left(\frac{1}{4^{r^{*}}n}\right)
\leq\frac{1}{3n}
\left(\sum_{i=1}^{n}|a_{i}|^{r}
\right)^{\frac{1}{r}}
\leq \frac{2}{3C_r}
N^{*-1}\left(\frac{1}{n}\right).
$$
Now we show (\ref{KombMittelLr73}).
Let $x\in\mathbb R^{n}$ with $\|x\|_{N^{*}}=1$
and $x_{1}\geq x_{2}\geq\cdots\geq x_{n}\geq 0$.
Furthermore, let $t\in\mathbb R^{n}$ s.t.
$x_{i}=N^{*-1}(t_{i})$. Let $i_{0}$ s.t.
$$
t_{1}\geq t_{2}\geq\cdots\geq t_{i_{0}}
\geq \frac{1}{4^{r^{*}}n}> t_{i_{0}+1}\geq\cdots\geq t_{n}.
$$
Then $\sum_{i=1}^{n}t_{i}=1$. We choose
$$
\tilde t=(t_{1},\dots,t_{i_{0}},0,\dots,0)
\hskip 20mm
{\tilde{\tilde{t}}}=(0,\dots,0,t_{i_{0}+1},\dots,t_{n})
$$
and
$$
\tilde x=(x_{1},\dots,x_{i_{0}},0,\dots,0)
\hskip 20mm
{\tilde{\tilde{x}}}=(0,\dots,0,x_{i_{0}+1},\dots,x_{n}).
$$
Then we have
$\|{\tilde{\tilde{x}}}\|_{N^{*}}\leq\frac{2}{3}$ and
\begin{equation}\label{KombMittelLr3}
\frac{1}{32\cdot 4^{r^{*}}}\| x\|_{N^{*}}
\leq \|\tilde x\|_{\bar N^{*}}
\leq 48\cdot 4^{r^{*}}\|\tilde x\|_{N^{*}}.
\end{equation}
We show this.
We estimate $\|{\tilde{\tilde{x}}}\|_{N^{*}}$.
We have
$$
\|{\tilde{\tilde{x}}}\|_{N^{*}}
=\inf\left\{\rho>0\left|\sum_{i=i_{0}+1}^{n}N^{*}
\left(\frac{N^{*-1}(t_{i})}{\rho}\right)\leq1\right.\right\}.
$$
By (\ref{KombMittelLr61})
$$
\sum_{i=i_{0}+1}^{n}N^{*}
\left(\frac{N^{*-1}(t_{i})}{\rho}\right)
\leq nN^{*}
\left(\frac{N^{*-1}(\frac{1}{4^{r^{*}}n})}{\rho}\right)
\leq nN^{*}
\left(\frac{2N^{*-1}(\frac{1}{n})}{3\rho}\right),
$$
and thus
$$
\|{\tilde{\tilde{x}}}\|_{N^{*}}\leq\frac{2}{3}.
$$
Therefore, $\|\tilde x\|_{N^{*}}\geq\frac{1}{3}$.
From (\ref{KombMittelLr2}) it follows
$$
\sum_{i=1}^{i_{0}}\bar N^{*}
\left(\frac{N^{*-1}(t_{i})}{\rho}\right)
\leq\sum_{i=1}^{i_{0}}\bar N^{*}
\left(\frac{16\cdot 4^{r^{*}}\bar N^{*-1}(t_{i})}{\rho}\right).
$$
Thus, we have
$$
\|\tilde x\|_{\bar N^{*}}\leq16\cdot 4^{r^{*}}.
$$
Using this and $\|\tilde x\|_{N^{*}}\geq\frac{1}{3}$ we obtain
$$
\|\tilde x\|_{\bar N^{*}}
\leq48\cdot 4^{r^{*}}\|\tilde x\|_{N^{*}}.
$$
Hence, the right hand side of inequality (\ref{KombMittelLr3}) is proved. Now we show the left hand side.
By (\ref{KombMittelLr74})
$$
\sum_{i=1}^{i_{0}}\bar N^{*}
\left(\frac{N^{*-1}(t_{i})}{\rho}\right)
\geq\sum_{i=1}^{i_{0}}\bar N^{*}
\left(\frac{\bar N^{*-1}(t_{i})}
{32\cdot 4^{r^{*}}\rho}\right).
$$
Thus
$\|\tilde x\|_{\bar N^{*}}\geq\frac{1}{32\cdot 4^{r^{*}}}$.
Using $\|x\|_{N^{*}}=1$, we obtain the left side of inequality (\ref{KombMittelLr3}), {\it i.e.},
$$
\|\tilde x\|_{\bar N^{*}}
\geq\frac{1}{32\cdot 4^{r^{*}}}\|x\|_{N^{*}}.
$$
The left hand side of inequality (\ref{KombMittelLr3}) implies the left hand side
of (\ref{KombMittelLr73}). The right hand side inequality of 
(\ref{KombMittelLr3}) implies 
$$
\|x\|_{\bar N^{*}}
\leq\|\tilde x\|_{\bar N^{*}}
+\|{\tilde{\tilde{x}}}\|_{\bar N^{*}}
\leq 48\cdot 4^{r^{*}}\|\tilde x\|_{N^{*}}
+\|{\tilde{\tilde{x}}}\|_{\bar N^{*}}.
$$
It is left to estimate the second summand. By
(\ref{KombMittelLr4})
$$
\sum_{i=i_{0}+1}^{n}\bar N^{*}
\left(\frac{N^{*-1}(t_{i})}{\rho}\right)
\leq n\bar N^{*}
\left(\frac{N^{*-1}(\frac{1}{n})}{\rho}\right)
\leq n\bar N^{*}
\left(\frac{16\bar N^{*-1}(\frac{1}{n})}{\rho}\right)
$$
Hence, $\|{\tilde{\tilde{x}}}\|_{\bar N^{*}}\leq16$
and 
$$
\|x\|_{\bar N^{*}}\leq(48\cdot 4^{r^{*}}+16)\|x\|_{N^{*}}.
$$
$\square$
\vskip 5mm

\begin{lem}\label{KombMittelLrLp}\label{THM Carsten 1} \label{THM Kwapien und Schuett}
Let $1\leq p<r<\infty$ and $a\in\mathbb R^{n}$ with
$a_{1}\geq a_{2}\geq\cdots a_{n}>0$. Then there exists an Orlicz function $N$
such that for the dual function $N^{*}$ and all $\ell=1,\dots,n$
\begin{eqnarray}\label{KombMittelLrLp1}
&& \tfrac{1}{2}
N^{*-1}\left(\frac{\ell}{n}\right)\\
 & \leq & 
C_r
\left\{
\left(\frac{\ell}{n}\right)^{\frac{1}{p^{*}}}
\left(\frac{1}{n}
\sum_{i=1}^{\ell}|a_{i}|^{p}\right)^{\frac{1}{p}}
+\left(\frac{\ell}{n}\right)^{\frac{1}{r^{*}}}
\left(\frac{1}{n}\sum_{i=\ell+1}^{n}|a_{i}|^{r}
\right)^{\frac{1}{r}}\right\}
 \leq  2^{-\frac{1}{p}}8
N^{*-1}\left(\frac{\ell}{n}\right)
\nonumber
\end{eqnarray}
holds, which is affine on the intervals 
$[\frac{\ell}{n},\frac{\ell+1}{n}]$, $\ell=0,1\dots,n-1$, and where $C_r=r^{1/r}(r^{*})^{1/r^*}$. For all such Orlicz functions and all
$x\in\mathbb R^{n}$ we have
$$
\alpha_{r,p}\|x\|_{N}
\leq\left(\frac{1}{n!}\sum_{\pi}
\left(\sum_{i=1}^{n}\left(|x(i)a_{\pi(i)}
\right|^{r}\right)^{\frac{p}{r}}
\right)^{\frac{1}{p}}
\leq \beta_{r,p}\|x\|_{N},
$$
where $\alpha_{r,p}$ and $\beta_{r,p}$ are constants, just depending on $r$ and $p$.
\end{lem}
\vskip 3mm

\noindent
{\bf Proof.}
For all $x\in\mathbb R^{n}$
$$
\frac{1}{n!}\sum_{\pi}
\|(x(i)a_{\pi(i)})_{i=1}^{n}\|_{r}^{p}
=\frac{1}{n!}\sum_{\pi}
\left(\sum_{i=1}^{n}\left(|x(i)|^{p}|a_{\pi(i)}|^{p}
\right)^{\frac{r}{p}}
\right)^{\frac{p}{r}}.
$$
Using Lemma \ref{KombMittelLr}, we get the existence of an Orlicz function
$M$ with
$$
a_{\frac{r}{p}}\|(|x(i)|^{p})_{i=1}^{n}\|_{M}
\leq
\frac{1}{n!}\sum_{\pi}
\left(\sum_{i=1}^{n}\left(|x(i)|^{p}|a_{\pi(i)}|^{p}
\right)^{\frac{r}{p}}
\right)^{\frac{p}{r}}
\leq b_{\frac{r}{p}}\|(|x(i)|^{p})_{i=1}^{n}\|_{M}
$$
and
\begin{eqnarray}\label{KombMittelLrLp1}
&&
M^{*-1}\left(\frac{\ell}{n}\right)   \\
&&\leq
C_r
\left\{\left(\frac{1}{n}
\sum_{i=1}^{\ell}|a_{i}|^{p}\right)
+\left(\frac{\ell}{n}\right)^{1-\frac{p}{r}}
\left(\frac{1}{n}\sum_{i=\ell+1}^{n}|a_{i}|^{r}
\right)^{\frac{p}{r}}\right\}
\leq 
4 M^{*-1}\left(\frac{\ell}{n}\right).
\nonumber
\end{eqnarray}
It follows
$$
(a_{\frac{r}{p}})^{\frac{1}{p}}
\|(|x(i)|^{p})_{i=1}^{n}\|_{M}^{\frac{1}{p}}
\leq
\left(\frac{1}{n!}\sum_{\pi}
\left(\sum_{i=1}^{n}\left(|x(i)|^{p}|a_{\pi(i)}|^{p}
\right)^{\frac{r}{p}}\right)^{\frac{p}{r}}
\right)^{\frac{1}{p}}
\leq
(b_{\frac{r}{p}})^{\frac{1}{p}}
\|(|x(i)|^{p})_{i=1}^{n}\|_{M}^{\frac{1}{p}}.
$$
Furthermore, we have
$$
\|(|x(i)|^{p})_{i=1}^{n}\|_{M}^{\frac{1}{p}}
=\|x\|_{M\circ t^{p}},
$$
since
$$
\|(|x(i)|^{p})_{i=1}^{n}\|_{M}^{\frac{1}{p}}
=\left\{\rho^{\frac{1}{p}}>0\left|
\sum_{i=1}^{n}M\left(\frac{|x(i)|^{p}}{\rho}\right)
\leq1\right.\right\}
=\left\{\eta>0\left|
\sum_{i=1}^{n}M\left(\left|\frac{x(i)}{\eta}
\right|^{p}\right)
\leq1\right.\right\}.
$$
We choose
$
N=M\circ t^{p}.
$
Then
$$
(a_{\frac{r}{p}})^{\frac{1}{p}}
\|x\|_{N}
\leq\left(\frac{1}{n!}\sum_{\pi}
\left(\sum_{i=1}^{n}\left(|x(i)|^{p}|a_{\pi(i)}|^{p}
\right)^{\frac{r}{p}}\right)^{\frac{p}{r}}
\right)^{\frac{1}{p}}
\leq (b_{\frac{r}{p}})^{\frac{1}{p}}\|x\|_{N}.
$$
Inequality (\ref{EQU M hoch -1 M stern hoch -1}) gives for all
$u\geq0$ 
$$
  u\leq M^{-1}(u)M^{*-1}(u)\leq 2u.
$$
Hence
$$
  N^{*-1}(t)\geq \frac{t}{N^{-1}(t)}
  =\frac{t}{(M^{-1}(t))^{\frac{1}{p}}}
  \geq2^{-\frac{1}{p}} t^{1-\frac{1}{p}}(M^{*-1}(t))^{\frac{1}{p}}
$$
$$
\tfrac{1}{2}N^{*-1}(t)\leq \frac{t}{N^{-1}(t)}
=\frac{t}{(M^{-1}(t))^{\frac{1}{p}}}
\leq t^{1-\frac{1}{p}}(M^{*-1}(t))^{\frac{1}{p}}
$$
Using (\ref{KombMittelLrLp1}), we get
\begin{eqnarray*}
&&\tfrac{1}{2}
N^{*-1}\left(\frac{\ell}{n}\right)    \\
&&\leq
C_r
\left\{
\left(\frac{\ell}{n}\right)^{\frac{1}{p^{*}}}
\left(\frac{1}{n}
\sum_{i=1}^{\ell}|a_{i}|^{p}\right)^{\frac{1}{p}}
+\left(\frac{\ell}{n}\right)^{\frac{1}{r^{*}}}
\left(\frac{1}{n}\sum_{i=\ell+1}^{n}|a_{i}|^{r}
\right)^{\frac{1}{r}}\right\}
\leq2^{-\frac{1}{p}}
N^{*-1}\left(\frac{\ell}{n}\right).
\end{eqnarray*}
$\square$

\vskip 3mm

The vector $(a_i)_{i=1}^n = \left(\frac{n}{i}\right)^{1/p}$ ($1<p<2$) generates the $\ell_p$-Norm, {\it i.e.}, we have for 
all $x\in\R^n$
\begin{equation}\label{GeneratingLp}
    c\norm{x}_p \leq \AveP \left( \sum_{i=1}^n \abs{x_ia_{\pi(i)}}^2\right)^{1/2} \leq C \norm{x}_p,
\end{equation}
where $c,C>0$ are absolute constants just depending on $p$.
This follows from Lemma \ref{KombMittelLrLp}.

\section{The Embedding of $\ell^n_M(\ell^n_N)$ into $L_1$}

To embed $\ell^n_M(\ell^n_r)$ into $L_1$, we have to extend the combinatorial expressions by another average over permutations. We use the following term
  \begin{equation}\label{EQU 3Fachmittel Produkt Raum}
    \AvePPP \left( \sum_{i,j=1}^n \abs{a_{ij}x_{\pi(i)}y_{\sigma(j)}z_{\eta(j)}}^2 \right)^{1/2}
  \end{equation}
which is, as we will show, under the appropriate choices of $x,y,z\in\R^n$ equivalent to  
  $$
    \norm{( \norm{(a_{ij})_{i=1}^n}_r )_{j=1}^n}_M.
  $$
Since (\ref{EQU 3Fachmittel Produkt Raum}) is equivalent to the $L_1$-Norm, we obtain the embedding into $L_1$. Using $z=((\frac{n}{j})^{1/p})_{j=1}^n$, $1<p<r<2$, we ``pass through'' an $\ell_{p}$ space to obtain the result. 

\begin{proposition} \label{THM Joscha u Carsten 2}
  Let $1<p<r<2$. Let $y\in\mathbb R^n\setminus\{0\}$ with $y_1 \geq \ldots \geq y_n \geq 0$,
  $(x_i)_{i=1}^n=((\frac{n}{i})^{1/r})_{i=1}^n$ and 
  $(z_j)_{j=1}^n=((\frac{n}{j})^{1/p})_{j=1}^n$. Then, for all matrices $a=(a_{ij})_{i,j=1}^n$,
\begin{eqnarray}\label{thmJoCa21}
a_{r,p} \|( \|(a_{ij})_{i=1}^n\|_r )_{j=1}^n\|_{M_y}
&\leq&
  \frac{1}{(n!)^{3}}\sum_{\pi,\sigma,\eta} \left( \sum_{i,j=1}^n |a_{ij}x_{\pi(i)}y_{\sigma(j)}z_{\eta(j)}|^2 \right)^{1/2} \nonumber \\
     &\leq& b_{r,p} \|( \|(a_{ij})_{i=1}^n\|_r )_{j=1}^n\|_{M_y},
\end{eqnarray}
where
$$
M_{y}\left(\frac{\ell}{n}\right)\sim
\frac{1}{n}\sum_{i=1}^{\ell}y_{i}
+\left(\frac{\ell}{n}\right)^{\frac{1}{p^{*}}}
\left(\frac{1}{n}\sum_{i=\ell+1}^{n}|y_{i}|^{p}\right)^{\frac{1}{p}}.
$$
  In particular, $\ell_{M_y}^n(\ell_r^n)$ is isomorphic to a subspace of $L_1$.  
\end{proposition} 

\noindent
{\bf Proof.}
  We start with the upper bound. By (\ref{GeneratingLp}) $z$ generates the $\ell_p$-norm. Thus
    \begin{equation}\label{thmJoCa2}
      \AvePPP \left( \sum_{i,j=1}^n \abs{a_{ij}x_{\pi(i)}y_{\sigma(j)}z_{\eta(j)}}^2 \right)^{1/2}
      \sim \AvePP \left( \sum_{j=1}^n y_{\sigma(j)}^p \left( \sum_{i=1}^n \abs{a_{ij}x_{\pi(i)}}^2 \right)^{p/2} \right)^{1/p}.
    \end{equation}
  By Jensen's inequality
    $$
      \AvePP \left( \sum_{j=1}^n y_{\sigma(j)}^p \left( \sum_{i=1}^n \abs{a_{ij}x_{\pi(i)}}^2 \right)^{p/2} \right)^{1/p}
      \leq \AvePs \left( \sum_{j=1}^n y_{\sigma(j)}^p \AveP \left( \sum_{i=1}^n \abs{a_{ij}x_{\pi(i)}}^2 \right)^{p/2} \right)^{1/p}.
    $$  
  By Lemma \ref{THM Carsten 1}, for all $j\leq n$
    $$
      \left( \AveP \left( \sum_{i=1}^n \abs{a_{ij}x_{\pi(i)}}^2 \right)^{p/2} \right)^{1/p} \sim \norm{(a_{ij})_{i=1}^n}_N,
    $$
  where 
    \begin{eqnarray*}
      N^{*-1}\left(\frac{\ell}{n}\right)
      & \sim& \left( \frac{\ell}{n} \right)^{1/p^*}\left(\frac{1}{n}\sum_{i=1}^l \abs{x_i}^p \right)^{1/p} 
      + \left( \frac{\ell}{n} \right)^{1/2}\left(\frac{1}{n}\sum_{i=\ell+1}^n \abs{x_i}^2\right)^{1/2}  \\
       & \sim& \left( \frac{\ell}{n} \right)^{1/p^*}\left(\frac{1}{n}\sum_{i=1}^{\ell} \left(\frac{n}{i}\right)^{p/r} \right)^{1/p}
      +   \left( \frac{\ell}{n} \right)^{1/2}\left(\frac{1}{n}\sum_{i=\ell+1}^n \left( \frac{n}{i} \right)^{2/r}\right)^{1/2} .
    \end{eqnarray*}
 Since $p<r<2$ 
     $$
     N^{*-1}(\tfrac{\ell}{n}) \sim (\tfrac{\ell}{n})^{1/r^*},
     $$ 
   which means that the N-norm is equivalent to the $\ell_{r}$-norm. Hence, we have shown the upper estimate of (\ref{thmJoCa21}),
  where $M_y$ is the Orlicz function as specified in Lemma \ref{THM Carsten 1}.\par
  For the lower bound, we obtain
    $$
      \AvePPP \left( \sum_{i,j=1}^n \abs{a_{ij}x_{\pi(i)}y_{\sigma(j)}z_{\eta(j)}}^2 \right)^{1/2}
      \sim \AvePP \left( \sum_{j=1}^n y_{\sigma(j)}^p \left( \sum_{i=1}^n \abs{a_{ij}x_{\pi(i)}}^2 \right)^{p/2} \right)^{1/p}.
    $$
  Now we use the triangle inequality and get
    $$
      \AvePP \left( \sum_{j=1}^n y_{\sigma(j)}^p \left( \sum_{i=1}^n \abs{a_{ij}x_{\pi(i)}}^2 \right)^{p/2} \right)^{1/p}
      \geq \AvePs \left( \sum_{j=1}^n y_{\sigma(j)}^p \abs{ \AveP \left( \sum_{i=1}^n \abs{a_{ij}x_{\pi(i)}}^2 \right)^{1/2}}^p \right)^{1/p}. 
    $$
  We know that for all $j\leq n$
    $$
      \AveP \left( \sum_{i=1}^n \abs{a_{ij}x_{\pi(i)}}^2 \right)^{1/2} \sim \norm{(a_{ij})_{i=1}^n}_r,
    $$
  since $(x_i)_{i=1}^n=((\frac{n}{i})^{1/r})_{i=1}^n$. Hence, by Lemma \ref{THM Carsten 1} we get the lower estimate of (\ref{thmJoCa21})
    $$
      \AvePPP \left( \sum_{i,j=1}^n \abs{a_{ij}x_{\pi(i)}y_{\sigma(j)}z_{\eta(j)}}^2 \right)^{1/2}
      \gtrsim \norm{( \norm{(a_{ij})_{i=1}^n}_r )_{j=1}^n}_{M_y}.
    $$
    Let $\epsilon$ and $\delta$ denote sequences of signs $\pm1$.
  Using (\ref{thmJoCa21}) and Khintchine's inequality one can easily show that
    $$
      \Psi_n : \ell^n_M(\ell_r^n) \to L_1^{n!^32^{2n}}, (a_{ij})_{i,j=1}^n \mapsto 
      \left( \sum_{i,j=1}^n a_{ij}x_{\pi(i)}y_{\sigma(j)}z_{\eta(j)}\e_i\delta_j \right)_{\pi,\sigma,\eta,\delta,\e} 
    $$
  embeds $\ell^n_M(\ell_r^n)$ into $L_1$.  
$\square$

\vskip 5mm

\begin{cor}
  Let $1<r<2$ and $1<p<r$. Furthermore, let $M$ be an $\alpha$-convex Orlicz function with $1<\alpha<p$. 
  Define $(x_i)_{i=1}^n=((\frac{n}{i})^{1/r})_{i=1}^n$, $(y_j)_{j=1}^n=(\frac{1}{M^{-1}(j/n)})_{j=1}^n$ and 
  $(z_j)_{j=1}^n=((\frac{n}{j})^{1/p})_{j=1}^n$. Then, for all matrices $a=(a_{ij})_{i,j=1}^n$,
    $$
      \AvePPP \left( \sum_{i,j=1}^n \abs{a_{ij}x_{\pi(i)}y_{\sigma(j)}z_{\eta(j)}}^2 \right)^{1/2}
      \sim \norm{( \norm{(a_{ij})_{i=1}^n}_r )_{j=1}^n}_{M}.
    $$
  In particular, $\ell_M^n(\ell_r^n)$ is isomorphic to a subspace of $L_1$.  
\end{cor}
\noindent
{\bf Proof.}
We apply Proposition \ref{THM Joscha u Carsten 2}.
We have to verify that the Orlicz function $M_{y}$ of 
Proposition \ref{THM Joscha u Carsten 2} is equivalent to
the Orlicz function $M$.
  We have for all $\ell\leq n$
    \begin{equation} \label{EQU Summation der ersten k �ber M-invers}
      \frac{1}{n} \sum_{i=1}^{\ell} \frac{1}{M^{-1}(\frac{i}{n})} \lesssim \frac{\ell}{n}\frac{1}{M^{-1}(\frac{\ell}{n})} 
      \stackrel{(\ref{EQU M hoch -1 M stern hoch -1})}{\sim} M^{*-1}\left( \frac{\ell}{n} \right)
    \end{equation}
  and
    \begin{equation}\label{EQU Summation der letzen �ber M-invers}
      \left( \frac{\ell}{n} \right)^{1/p^*} \left(\frac{1}{n} \sum_{i=\ell+1}^n \abs{\frac{1}{M^{-1}(\frac{i}{n})}}^p \right)^{1/p}
      \stackrel{(\ref{EQU M hoch -1 M stern hoch -1})}{\lesssim} M^{*-1}\left( \frac{\ell}{n} \right),
    \end{equation}  
  since $M$ is $\alpha$-convex and therefore $(M^{-1})^{\alpha}$ concave, {\it i.e.}, for all $l\leq n$
    $$
      M_y^{*-1}\left(\frac{\ell}{n}\right) \sim \frac{1}{n}\sum_{i=1}^{\ell} y_i  
    + \left( \frac{\ell}{n} \right)^{1/p^*}\left(\frac{1}{n}\sum_{i=\ell+1}^n \abs{y_i}^p\right)^{1/p}\lesssim M^{-1}\left(\frac{\ell}{n}\right).
    $$
  The lower bound is trivial, since $M^{-1}$ is an increasing function. 
$\square$

\vskip 5mm

We will now prove Theorem \ref{MatrixEmbed}.

\vskip 5mm

\noindent
{\bf Proof of Theorem \ref{MatrixEmbed}.}
It is enough to show the case $N(t)=t^{r}$.
Indeed, by \cite{key-1} $\ell_{N}$ is a subspace of 
$L_{r}$ if and only if $\frac{N(t)}{t^{r}}$ increasing and 
$\frac{N(t)}{t^{2}}$ decreasing.
\par
We apply Proposition \ref{THM Joscha u Carsten 2}.
We choose $y_{i}$, $i=1,\dots,n$, such that
$$
M^{*-1}\left(\frac{\ell}{n}\right)=
\frac{1}{n}\sum_{i=1}^{\ell}y_{i}.
$$
We show that $M^{*}$ and $M^{*-1}_{y}$ of
Proposition \ref{THM Joscha u Carsten 2} are equivalent.
It follows that for all $\ell$
$$
M^{*-1}\left(\frac{\ell}{n}\right)
\geq\frac{\ell}{n}y_{\ell}
$$
By
$$
t\leq M^{-1}(t)M^{*-1}(t)\leq 2t,
$$
we get
$$
y_{\ell}
\leq\frac{M^{*-1}\left(\frac{\ell}{n}\right)}
{\frac{\ell}{n}}
\leq\frac{2}{M^{-1}\left(\frac{\ell}{n}\right)}.
$$
Therefore,
\begin{eqnarray*}
\frac{1}{n}\sum_{i=\ell+1}^{n}|y_{i}|^{p}
\leq2^{p}\frac{1}{n}\sum_{i=\ell+1}^{n}
\frac{1}{|M^{-1}\left(\frac{i}{n}\right)|^{p}}
=2^{p}\frac{1}{n}\sum_{i=\ell+1}^{n}
\frac{|\frac{i}{n}|^{\frac{p}{r}}}
{|M^{-1}\left(\frac{i}{n}\right)|^{p}} 
|\frac{n}{i}|^{\frac{p}{r}} . 
\end{eqnarray*}
Since $\frac{M(t)}{t^{r}}$ is decreasing, 
$\frac{s}{|M^{-1}(s)|^{r}}$ is decreasing. Therefore,
since $r<p$
$$
\left|\frac{t}{|M^{-1}(t)|^{r}}\right|^{\frac{p}{r}}
=\frac{t^{\frac{p}{r}}}{|M^{-1}(t)|^{p}}
$$
is also decreasing. Thus
\begin{eqnarray*}
\frac{1}{n} 
\sum_{i=\ell+1}^{n}|y_{i}|^{p}
&\leq&2^{p}
\frac{|\frac{\ell}{n}|^{\frac{p}{r}}}
{|M^{-1}\left(\frac{\ell}{n}\right)|^{p}}\frac{1}{n}\sum_{i=\ell+1}^{n}
|\frac{n}{i}|^{\frac{p}{r}} 
\leq2^{p}
\frac{|\ell|^{\frac{p}{r}}}
{|M^{-1}\left(\frac{\ell}{n}\right)|^{p}}\frac{1}{n}\sum_{i=\ell+1}^{n}
i^{-\frac{p}{r}}  \\
&\sim&2^{p}
\frac{|\ell|^{\frac{p}{r}}}
{|M^{-1}\left(\frac{\ell}{n}\right)|^{p}}\frac{1}{n}
\ell^{1-\frac{p}{r}}
=2^{p}\frac{\ell}{n}
\frac{1}
{|M^{-1}\left(\frac{\ell}{n}\right)|^{p}}.
\end{eqnarray*}
Altogether,
\begin{eqnarray*}
M^{*-1}\left(\frac{\ell}{n}\right)
&\leq&
\frac{1}{n}\sum_{i=1}^{\ell}y_{i}
+\left(\frac{\ell}{n}\right)^{\frac{1}{p^{*}}}
\left(\frac{1}{n}\sum_{i=\ell+1}^{n}|y_{i}|^{p}\right)^{\frac{1}{p}}  \\
&\leq&M^{*-1}\left(\frac{\ell}{n}\right)
+\left(\frac{\ell}{n}\right)^{\frac{1}{p^{*}}}
\left(2^{p}\frac{\ell}{n}
\frac{1}
{|M^{-1}\left(\frac{\ell}{n}\right)|^{p}}\right)^{\frac{1}{p}}   \\
&=&M^{*-1}\left(\frac{\ell}{n}\right)
+2\frac{\ell}{n}
\frac{1}
{|M^{-1}\left(\frac{\ell}{n}\right)|^{}}
\leq 3M^{*-1}\left(\frac{\ell}{n}\right).
\end{eqnarray*}
$\square$

$$
  \begin{array}{ll}
    \hbox{J. Prochno} & \hbox{C. Sch\"utt} \\
    \hbox{Christian-Albrechts-Univerist\"at zu Kiel} & \hbox{Christian-Albrechts-Univerist\"at zu Kiel} \\
    \hbox{Ludewig-Meyn-Str. 4} & \hbox{Ludewig-Meyn-Str. 4} \\
    \hbox{24098 Kiel} & \hbox{24098 Kiel} \\
    \hbox{prochno@math.uni-kiel.de} & \hbox{schuett@math.uni-kiel.de} \\
  \end{array}
$$

\end{document}